# Identify the Optimal Locations of Dedicated Bus Lanes to Improve the Stability of Bus Line

He Sheng-Xue[1*]

1. Business School, University of Shanghai for Science and Technology, Shanghai 200093, China

2. School of Business, State University of New York at Oswego, 7060 State Route 104, NY, 13126-3599, USA

**Abstract**: Bus bunching and unevenly dispersed buses along a bus line lead to a low service level and deteriorate the operational stability of the bus line. Speed adjustment as a control means has been proposed to solve the above problem. But due to the difficulty of adjusting speed in the mixed traffic and the discontinuously distributed dedicated bus lanes (DBLs) in a bus line, the existing methods which coordinate or adjust bus speeds in all road segments encounter serious obstacles in practice.

To overcome the above problem, we first take into account the influence of the deployment of the DBLs along a bus line on the performance of the speed adjustment strategy. With the funding limit and the maximally allowable influence on other road traffic as main constraints, a nonlinear mixed integer stochastic mathematic model is formulated aiming at improving the whole stability of the bus line operation. The model can be effectively solved by a newly designed branch and bound algorithm. The optimal regulating speed in a DBL is determined with a process of looking multiple critical time point ahead.

Numerical experiment verified the effectiveness and efficiency of our new method. Three findings were obtained as follows. To determine the optimal regulating speed, the number of the critical time points to be looked ahead should not be too big. Increasing the number of DBLs improves the stability of a bus line. The more evenly DBLs are distributed along the bus line, the greater the stability of the bus line will be.

**Keywords**: traffic engineering; modeling; transport planning; bus bunching; dedicated bus lane; branch and bound algorithm

## 1. Introduction

Dedicated Bus Lanes (DBL) are widely accepted as an effective way to improve the public transit service. On the one hand, we see the average cruising speeds of buses are increased after the implementation of DBLs along a bus line (Zyryanov and Mironchuk, 2012; Ben-Dor, et al. 2018). On the other hand, the service level with respect to the operational instability and reliability has not shown an expected improvement in many cases (Arason and Vedagiri, 2010). Sometimes passengers may experience a more congested ridding process and need to wait longer for a bus after the DBLs are adopted in a bus line. To uncover the underlying reason and solve the above problem, further researches are required.

Two decisions will together influence the final function of DBLs. The first one is where to install DBLs under some necessary constraints. The constraints include the funding limit and the maximally allowable influence on the other traffic. The second one is that after the locations of DBLs have been determined, how to make use of these special infrastructures to improve the operational stability and reliability of a bus line. Since buses have the priority of road right in DBLs, it is usually feasible and convenient to adjust the speeds of buses in a proper range in a DBL (Arason and Vedagiri, 2010; Zyryanov and Mironchuk, 2012; Ben-Dor, et al. 2018). So making use of DBLs could first take into account the speed adjustment. We will focus on these two issues in this paper.

---

[*] Tel.: +086-18019495040.
E-mail address: lovellhe@126.com

Researchers have made many efforts to improve the operational reliability and stability of a bus line. The proposed methods include stop-skipping strategies (Suh et al. 2002; Fu et al. 2003; Sun and Hickman 2005; Cortés et al. 2010; Liu et al. 2013; Zhang et al. 2019), limited-boarding strategies (Osuna and Newell 1972; Newell 1974; Barnett 1974; Delgado et al., 2009, 2012), embedding-slack strategies (Daganzo 1997a, b; Daganzo 2009a; Xuan et al. 2011), static and dynamic holding strategies (Hickman 2001; Eberlein et al. 2001, Sun and Hickman 2008; Puong and Wilson 2008; Daganzo 2009a; Xuan et al. 2011; Bartholdi and Eisenstein 2012; Delgado et al. 2012; He 2015; Argote-Cabanero et al. 2015), coordinated speed adjustment strategies (Daganzo 2009b; Daganzo and Pilachowski 2009, 2011), transit signal priority strategies (Liu et al. 2003; Ling and Shalaby 2004; Estrada et al. 2016) and bus substitution strategy (Petit et al. 2018).

Among the existing literature, only a few studies focused on the speed adjustment. Researchers tried to resist bus bunching and to improve the operational instability of a bus line with the speed adjustment (Daganzo 2009b; Daganzo and Pilachowski 2009, 2011; He 2015). By dynamically coordinating the speeds of successive cruising buses, Daganzo and Pilachowski (2009, 2011) presented a speed adjustment strategy which can effectively cope with the big fluctuation or interference on a bus line. But as they mentioned, to implement their strategy in practice may encounter some challenges. As we all know if a bus travels among the congested traffic, it is hard for a bus driver to adjust the bus speed at will. Even if the DBLs are considered, it is still very difficult to coordinate the speeds of buses which are cruising along a bus line with discontinuously dispersed DBLs. In many cases, we can only adjust the speeds of buses in DBLs. The speed of buses in the common road segments cannot be adjusted. He (2015) suggests adjusting the bus speeds on all general road segments only at some critical time points. Though the method in He (2015) does not need to coordinate the speeds of the successive buses dynamically, its implementation is limited by the same constraint mentioned above to Daganzo and Pilachowski (2009, 2011).

The existing studies about the deployment of dedicated bus lane focused on the influences of the installation of DBLs on the travel times, pollutant emission and operational cost (Sun and Wu 2017). Most of them dealt with the deployment problem from a network viewpoint and used the traffic assignment model as the basis of analyses (Miandoabchi et al. 2012; Khoo et al. 2014; Sun and Wu 2017). The existing studies seldom dealt with the connection of DBLs' deployment with the improvement of stability of a bus line. The above lack may lead to an unexpected negative performance to a bus line. A study focusing on the above connection has the potential to further improve the stability of a bus line and avoid some irritating side effects of installing DBLs.

Selecting the locations of DBLs and making use of DBLs entwine together. To make any one decision above being effective and practical, the other decision must be considered carefully. As far as the authors know, most of the existing researches didn't address them together. The lack of a systematical study about their integrated influence makes it hard to realize their potential at improving the system stability. We need to know where to install the DBLs and then how to improve the stability of a bus line by adjusting the cruising speeds of buses in the installed DBLs.

To answer the above questions, we first present a mixed integer nonlinear stochastic model to formulate our goal and the decision setting. The operational instability can be measured properly by our new defined index. Our goal is to minimize the whole instability of a bus line during a given observation period. The funding limit and the maximally allowable influence on the other traffic will be considered as the main constraints in the model.

To resolve the obtained model, we divide our task into two subtasks. The first subtask is to design an effective strategy to minimize the system instability during the observation period with the given locations of DBLs. After estimating and comparing the influences of all alternatives of the speed variable on the stability of the bus line with a multiple Critical Time Point (CTP) looking ahead mechanism, a proper regulating speed is chosen. The above mechanism uses a detailed simulation system to describe the complex bus line components, such as the random arrivals of passengers at stops and the signalized intersections. The above strategy can cope with the discontinuously dispersed DBLs along a bus line since this strategy requires no compulsive coordination of speeds among the successive cruising buses.

The second subtask is to determine the optimal locations for installing DBLs with the help of the results of the first subtask. With the strategy of solving the first subtask, we can choose the best locations by comparing the expected optimal performances of all the alternative combinations of feasible locations. The essence of this problem is to choose a combination of elements in a given set. But in view of the time consuming process, this is not easy if the size of the set is big. To overcome the above obstacle, we design a specific branch and bound algorithm. By gradually reducing the number of the locations where the installation of DBL is allowed, a searching tree will be generated. This algorithm has the power to go through all the combinations of the possible locations, but can wisely avoid unnecessary searching.

By putting the two subtasks together, we can resolve our problem effectively. A bus line simulated in detail will be used to verify our approach. Analyzing the numerical examples will uncover some unknown insights about the location selection and the speed adjustment.

The remainder of this paper is organized as follows. Section 2 formulates the location selection model about DBLs' deployment. Section 3 introduces the branch and bound algorithm used to determine the optimal locations of DBLs. Section 4 presents the method of looking multiple critical time points ahead to choose a proper regulating speed for a bus cruising in a DBL. Section 5 uses examples to verify the effectiveness and efficiency of our model and method. Section 6 sums up this study and points out the possible directions for future researches.

## 2. Location Selection Model of DBLs

### 2.1. Instability Indices and Objective

*2.1.1. Control action*

# Position of figure 1

Fig. 1 The partial trajectories of two buses

In Figure 1, the partial trajectories of two buses $b_1$ and $b_2$ are depicted. Assume that there are some dedicated bus lanes connecting bus stop $e_1$ with $e_2$. When bus $b_1$ leaves from $e_1$ at time $t_3$, its cruising speed can be changed from the average expected speed corresponding to the solid slant line to the higher speed corresponding to the dot-dash slant line. The speed adjustment operation as above is the only control action to be considered in this paper. Note that we will also consider slowing down a bus in a DBL if necessary. From Figure 1, we can see the headway of bus $b_1$ to its preceding bus $b_2$ at time $t_3$ is $h_1(t_3)$ after the speed of $b_1$ is adjusted. If the speed of $b_1$ remains the average expected speed, the headway should be $h_1(t_3) + \Delta h$.

Let $l$ denote a dedicated bus lane (DBL). Assume that the average cruising speed of buses in $l$ is $\bar{v}_l$. When a bus travels in this DBL, its speed may be adjusted to stabilize the system. Its new speed is denoted by $\bar{v}_l + \Delta v_l$ where $\Delta v_l$ is the regulating speed. To make the control practical, we assume that $\Delta v_l$ can only be chosen from a given set $A_l$ which includes a finite number of elements. $A_l$ is called the action set with respect to $l$. For example, $A_l$ may be composed of -15, -10, -5, 0, 5, 10, and 15. The measurement unit of $\Delta v_l$ is the same as speed, i.e. km/hr in this paper. We will use $a$ or $a_l$ to stand for a feasible value of $\Delta v_l$ in $A_l$. $a$ is called a control action. If we view the set {0} with only one element 0 as the action set of all the common non-DBL road segments, all the road segments with or without DBL can be dealt with in a similar way in the following expressions.

A bus line segment is the part of a bus line connecting two neighboring bus stops. A bus line segment generally consists of several road segments which are the part of a bus line with two neighboring bus stops or intersections as its ends. In the real world, the road segments with DBLs may not connect with each other. To simplify the following analysis, we assume that the regulating speeds on all the dedicated bus lanes included in a bus line segment will be the same. In other words, the control action $a$ with respect to a bus line segment is the same.

*2.1.2. Measurement of Instability*

To survey the stability of a dynamically changing bus line system during a given time period, we will choose some Critical Time Points (CTP) as the surveillance time points and note down the performance of the system at these points. To plot the trajectory of a bus, we need its arrival and departure times with respect to the bus stops which have been passed by. Since the regulating speed needs to be determined before the bus enters into a bus line segment with dedicated bus lanes, we will use the departure times as our surveillance time points. Obviously, the speed adjustment decision will influence the following headways between some buses cruising along the bus line. Even if there is no speed adjustment when a bus leaves a stop, we will still calculate the headways between buses at this time point. These obtained headways will be used later to evaluate the stability performance of the bus line. Note that during the observation period, it may appear that several buses depart from their dwelling stops at the same time. In this case, the critical time point is the same for these buses. Since we may need to adjust their speeds at the same time in this case, the order to deal with them will be arbitrary. Since this phenomenon is very rare, its impact on the final result is negligible.

To know whether a bus line system runs stable at a CTP or not, we can check if the instantaneous headways between the successively cruising buses are evenly distributed. Use $h_b(t)$ to denote the headway between bus $b$ and its nearest preceding bus at time $t$. The time headway is the time required for a bus traveling from its current position to the current position of its nearest preceding bus. The Dynamic Circle Headway (DCH) with respect to time $t$ is defined as follows:

$$H(t) = \sum_{b \in B} h_b(t) / n_b, \qquad (1)$$

where $B$ is the set of all buses and $n_b$ is the total number of buses. DCH is an average value at the given time $t$. Its corresponding standard deviation is given below:

$$\sigma_H(t) = \sqrt{[\sum_{b \in B}(h_b(t) - H(t))^2]/n_b} \quad . \qquad (2)$$

$\sigma_H(t)$ indicates the deviation of headways from the average value. So we can use $\sigma_H(t)$ to estimate the stability of a bus line at time $t$.

To estimate the whole stability of a bus line system during the given observation period, we need to define some new index over the observation period with $\sigma_H(t)$. Let $\bar{T}$ denote the set of all critical time points corresponding to the departure times of buses from bus stops. The whole stability index is defined as follows:

$$\bar{c}_H = (\sum_{t \in \bar{T}} \sigma_H(t))/n_{\bar{T}}, \qquad (3)$$

where $n_{\bar{T}}$ is the cardinality of $\bar{T}$.

$\bar{c}_H$ stands for the deviation of headways from DCH over the whole observation period. From now on, $\bar{c}_H$ will be referred to as the First Stability Index (**FSI**). Since $\bar{c}_H$ is defined as average values over the cruising buses and the set $\bar{T}$ of critical time points, it intrinsically possesses the potential to indicate the operational stability of any bus line system as a whole.

*2.1.3. Problem Description and the Objective*

In the real world, it is very rare to install DBLs for all the road segments of a bus line. Some road segments may be unfit for installing DBLs. For example, a road segment has only one unidirectional lane. Even if there are several lanes in a road segment, the heavy traffic may prohibit the installation of any DBL. If there are many road segments which are fit for the installation, the final influence of any deployed scheme needs to be evaluated after a careful comparison among different plans.

To make the following formulation simple, we assume that an installing or not decision is always corresponding to a bus line segment. Obviously, a bus line segment may include several intersections and road segments. Sometimes all the road segments included in a bus line segment are fit for installing DBLs. Sometimes they are all unfit. No matter what is the real situation, we group the decisions of installing or not for all the road segments into one decision. The above simplicity will not change the main theoretical analyses later. This simplification assumption can be easily fixed by extending the model properly based on road segments.

Let $I$ denote the ordered set of all feasible locations where DBLs can be installed. Location here means the bus line segments in this paper. One typical element of $I$ is indicated by $i$. The decision to install DBLs or not in the bus line segment $i$ is expressed by variable $x_i \in \{0,1\}$. If $x_i$ equals 1, the decision is to install the related DBLs; or else, the decision is not to install any DBL. Group all the decision variables $x_i$, $\forall i$ into a column vector $x$ which corresponds to the elements' order of set $I$.

Use $T$ to denote the observation period. Let $a_m$ be the control action or the regulating speed at the critical time point $t_m$. Use $a(x,T)$ or $a(x)$ to denote the vector which consists of all the control actions, such as $a_m$, during $T$ on a bus line with a given value of $x$. Note that the number of elements in $a(x)$ is generally indeterminate until the end of the observation.

In eq. (3), we define the stability index $\bar{c}_H$. The value of $\bar{c}_H$ depends on the given pattern $x$ of DBLs' locations and the control actions $a(x)$. We use $\bar{c}_H(x,a)$ to denote the above relationships. Minimizing $\bar{c}_H(x,a)$ is equivalent to stabilizing a bus line with $x$ and $a(x)$ as the decision variables. For a dynamically evolving bus line system, a speed adjustment decision is made at the critical time point $t_m$ based on the system state $s_m$ and the given value of $x$. We call the rules to make the above decision a policy $\pi$. The decision-making process is denoted by function $a_m(x) = A^\pi(s_m, x)$. $s_m$ will be defined and explained in detail later in subsection 4.1. $s_m$ consists of the main features of a bus line at the time $t_m$. To make the subsequent expressions concise, we group all the relations $a_m(x) = A^\pi(s_m, x)$, $\forall t_m$ into a mapping $a(x) = A^\pi(s, x)$.

With the above notations, we can define our optimization objective as follows:

$$\min \mathbb{E}[\bar{c}_H(x, A^\pi(s, x))], \qquad (4)$$

where the operator "$\mathbb{E}$" means to obtain the expected value of a random variable. The randomness of the above objective comes from the various random and uncertain factors of a bus line, such as the random arrivals of passengers at stops and the uncertain influence of other road traffic on the cruising buses.

**2.2. Constraints**

It is not free to install a dedicated bus lane. There are some constraints required to be taken into account. In this paper, we mainly consider two types of constraints. The first one is about the influence of the installation of DBLs on the other road traffic. Corresponding to the bus line segment $i$ in the feasible ordered set $I$, we use $F_i^1$ to denote the influence of installing DBLs on the other road traffic. And we assume that $F_i^1$ can be estimated and will be given as a constant in our model. All the $F_i^1$, $\forall i \in I$ constitute a column vector $F^1$. The upper limit of the total influence on the other road traffic is denoted by $\hat{F}^1$.

The second constraint is about the limitation of funds. Let $F_i^2$ be the necessary investment required to install DBLs on the bus line segment $i \in I$. All the necessary investment $F_i^2$, $\forall i \in I$ constitute a column vector $F^2$. Let $\hat{F}^2$ be the total feasible funds for the installation of DBLs.

To make the expression concise, we construct a matrix $F = \{F^1, F^2\}^\perp$ and a column vector $\hat{F} = \{\hat{F}^1, \hat{F}^2\}^\perp$ where the superscript "$\perp$" stands for the transposition operation. With the above notations, the constraints can be given as follows:

$$Fx \leq \hat{F}. \tag{5}$$

Obviously, except the constraints (5), the elements of variable $x$ should also be limited to the binary set $\{0,1\}$.

**2.3. Basic Thought of Resolving the Model**

The model formulated above with the function (4) as the objective function is a mixed integer nonlinear stochastic programming model. Nowadays there exists no commonly accepted method to solve such kind of models efficiently and effectively. So we need to design a specific algorithm for our model. The basic thought to solve the above problem is explained as follows.

We will realize our goal in two steps. At first, we assume that the locations to install DBLs have been given. Our task now is to make speed adjustments by some effective strategy to minimize the stability index. In other words, we aim to minimize the value of $\bar{c}_H(x, A^\pi(s, x))$ by a policy $\pi$ with $x$ given. Here a proper policy $\pi$ will be the key to the minimization. We will deal with this issue later in Section 4. Since a bus line system is a stochastic system, the expected value $\mathbb{E}[\overline{c}_H(x, A^\pi(s,x))]$ will be estimated by an average value over a given number of the calculations of $\overline{c}_H(x, A^\pi(s,x))$.

The second step to realize our goal is to determine the optimal value of $x$. From the first step, we can obtain the expected value of $\mathbb{E}[\overline{c}_H(x, A^\pi(s,x))]$ for a given $x$. By comparing the values of $\mathbb{E}[\overline{c}_H(x, A^\pi(s,x))]$, the performances of two different values of $x$ can be evaluated. The smaller expected value means a more stable system. If the size of the set of the feasible values of $x$ is relatively small, we can get the optimal one by enumeration method. Unfortunately, when the number of elements in $x$ is large, the set of the feasible values of $x$ will be unthinkably large. For example, if there are 15 elements in the set, the possible values of $x$ will be 32,768 (that is $2^{15}$). To calculate the expected values of $\mathbb{E}[\overline{c}_H(x, A^\pi(s,x))]$ for all the feasible 32,768 values of $x$ will be time-consuming. Especially when the computational time required to obtain the value $\overline{c}_H(x, A^\pi(s,x))$ for a given $x$ is long, the comparison among all the feasible $x$ will be impractical. To overcome the above obstacle, Section 3 will present a specific branch and bound algorithm to considerably reduce the searching time.

## 3. Branch and Bound Algorithm

### 3.1. Basic Thought

The underlying idea of our branch and bound algorithm is to generate a searching tree from a root node. The root node stands for the choice that all the feasible locations are chosen to install DBLs. The root node corresponds to $x_i = 1$, $\forall i \in I$. A branch will spring out from a node by removing one chosen location from the set of the node's chosen locations. With the generation of a new branch, a new node as the head node of the newly generated branch is also generated. The reason why we start from a root node where all feasible locations are chosen and then gradually remove some chosen locations is that by doing so the objective value will be gradually increased. Note that by removing a chosen location, the capacity of speed adjustment will be naturally reduced.

With the shrinking of the set of the chosen locations, the constraints of (5) may be satisfied at last. When the constraints are satisfied, a feasible solution is obtained. A feasible solution supplies an upper limit to the subsequent searching. The more feasible solutions are found, the smaller the upper limit may be.

When does the searching or branching operation need to retrace the parent nodes? One situation is when a feasible solution is found. It is easy to see that to reduce an element from the set of the chosen locations of a feasible solution will further increase the objective value. The other situation is when the objective value corresponding to the current node is greater than the known upper limit of the known feasible solutions. In this situation, any further branching from current node will not further improve the obtained result.

### 3.2. Introduction of Related Notations

To specify the branch and bound algorithm, we need to clarify some notations related to the nodes $\mathbf{N}$ and the branches $\mathbf{A}$ of the searching tree $Tree(\mathbf{N}, \mathbf{A})$. From any node $n \in \mathbf{N}$ in the searching tree $Tree(\mathbf{N}, \mathbf{A})$ to the root node $n_0$, there is a unique path consisting of the removed locations. These removed locations are indicated by the traversed branches. By removing these locations, the left locations will constitute the set of the chosen locations for node $n$.

For any node $n$, there is an ordered list of alternative locations $\mathbf{L}_n$ corresponding to the above-mentioned set of the chosen locations. This list includes all the locations where to install DBLs is allowed at the current node. With the given ordered list $\mathbf{L}_n$, the corresponding objective value can be calculated by the strategy to be presented in Section 4.

Corresponding to $\mathbf{L}_n$, there is an ordered list $\bar{\mathbf{L}}_n$ which consists of all the current feasible alternative locations. But with the branching operation, the size of $\bar{\mathbf{L}}_n$ will be reduced gradually. A branch emitting from node $n$ will be distinguished by the location removed from the ordered list $\bar{\mathbf{L}}_n$. With the increase of branches emitting from a node, the number of its current feasible alternative locations in $\bar{\mathbf{L}}_n$ will decrease.

If $\mathbf{L}_n$ is a feasible solution to our model, node $n$ should be labeled by "Feasible". If all the possible branches emitting from a node $n$ have been generated, the node should be labeled by "Checked" and at this time $\bar{\mathbf{L}}_n$ becomes an empty set. Obviously, except the root node, any other node has a parent node connecting to itself by a branch.

A branch labeled by $g \in \mathbf{A}$ corresponds to a removed location $i \in I$ that is a bus line segment. So we use $i_g$ to denote the removed location associated with the branch $g$. Any branch connects a parent node to one of its descendant nodes.

### 3.3. The Branch and Bound Algorithm

Assume that there exists at least one feasible solution to our location problem.

**Step 1**: Initialization. Let $n_0$ be the root node such that $\mathbf{L}_{n_0} := I$ and $\bar{\mathbf{L}}_{n_0} = I$. Calculate the objective value $CH(\mathbf{L}_{n_0})$ of $\mathbf{L}_{n_0}$. And let $n_0$ be the current node $n$. Let the currently feasible solution $\tilde{L}$ be an empty set $\Phi$. Use $\Gamma$ to denote the upper limit for the following searching. Note that the initial value of $\Gamma$ should be initialized by a sufficiently large positive number.

**Step 2**: Bounding operation and termination judgment.

**Step 2.1**: If $CH(\mathbf{L}_n) \geq \Gamma$, remain the current values of $\tilde{L}$ and $\Gamma$ and let the parent node of $n$ be the current node. With the renewed current node $n$, restart Step 2.

**Step 2.2**: If the constraints of (5) are satisfied by $\mathbf{L}_n$ and $CH(\mathbf{L}_n) < \Gamma$, let $\mathbf{L}_n$ be the currently feasible solution $\tilde{L}$ and assign $CH(\mathbf{L}_n)$ to the upper limit $\Gamma$. If $n$ is the root note $n_0$, go to Step 4; or else, let the parent node of $n$ be the current node. If the current node is renewed, restart Step 2.

**Step 2.3**: If the constraints of (5) are not satisfied by $\mathbf{L}_n$ and $CH(\mathbf{L}_n) < \Gamma$, go to Step 3.

**Step 3:** Branching.

**Step 3.1:** If $\bar{\mathbf{L}}_n$ equals $\Phi$ and $n \neq n_0$, let the parent node of $n$ be the current node and go to Step 2. If $\bar{\mathbf{L}}_n$ equals $\Phi$ and $n = n_0$, go to Step 4.

**Step 3.2:** If $\bar{\mathbf{L}}_n$ is not equal to $\Phi$, choose the first element $i_n$ of $\bar{\mathbf{L}}_n$ as the location to be removed. Update $\bar{\mathbf{L}}_n$ by $\bar{\mathbf{L}}_n / \{i_n\}$. Then generate a descendent node $\eta$ of $n$. The node $\eta$ has $\mathbf{L}_n / \{i_n\}$ and the updated $\bar{\mathbf{L}}_n$ as its ordered sets $\mathbf{L}_\eta$ and $\bar{\mathbf{L}}_\eta$, respectively. A new branch $i_n$ connecting $\eta$ to $n$ is also produced. Calculate the objective value $CH(\mathbf{L}_\eta)$ of $\eta$. Let $n := \eta$, go to Step 2.

**Step 4**: Output the optimal locations $\tilde{L}$ and its corresponding objective value $\Gamma$. Set $x_i = 1$, $\forall i \in \tilde{L}$ and $x_i = 0$, $\forall i \notin \tilde{L}$. These renewed $x_i$, $\forall i \in I$ constitute the optimal choice $x^*$.

### 3.4. Completeness of the Algorithm

**Position of figure 2**

Fig. 2 A searching tree with 4 locations to be chosen at the beginning

**Proposition 1:** The algorithm given above can get the optimal solution of the location problem by checking all possible solutions.

**Proof**: As we know, a possible solution is a combination of set $I$. If we can verify the following two statements, this proposition will be proved. The first statement is that any node of the searching tree corresponds to a unique combination of set $I$. The second statement is that the searching process can generate a tree which includes all the possible nodes. A node here corresponds to a possible combination of $I$. Assume that the ordered set $I$ is $\{1, 2, \ldots, i, \ldots, l_I\}$. Any one of its combinations can be expressed by an ordered list $\{j_1, j_2, \ldots, j_K\}$ such that $1 \leq j_1 < j_2 < \ldots < j_K$ and $K \leq l_I$. All the elements of $\{j_1, j_2, \ldots, j_K\}$ come from the set $I$. Now we can show that $\{j_1, j_2, \ldots, j_K\}$ can be related to a unique path from the root node to a possible node of the searching tree. At first, $j_1$ will correspond to a branch emitted from the root node. All the elements of $\{i \in I | i < j_1\}$ should have been considered by our algorithm. So the ordered list $\bar{\mathbf{L}}_{n(j_1)}$ of the alternative locations connecting $j_1$ with $j_2$ will be $\{j_1+1, j_1+2, \ldots, l_I\}$. Then, $j_2$ will correspond to a branch emitted from the node $n(j_1)$. All the elements of $\{i \in I | i < j_2\}$ should have been considered by our algorithm. So the ordered list $\bar{\mathbf{L}}_{n(j_2)}$ of alternative locations connecting $j_2$ with $j_3$ will be $\{j_2+1, j_2+2, \ldots, l_I\}$. The above process can be continued until generating a specified path. The ordered branches included in this path are indicated by $\{j_1, j_2, \ldots, j_K\}$. The above analysis has clarified the one-to-one corresponding relations between the combination of set $I$ and the nodes of the searching tree. The statement is proved. □

The tree plotted in Figure 2 can be used to facilitate the understanding of the above proof. Note that the tree in Figure 2 is generated from left to right. $\{1, 2, 3, 4\}$ is the set $I$ in Figure 2. The series in $\{\ldots\}$ and $(\ldots)$ beside the nodes of Figure 2 are used to give $\mathbf{L}_n$ and the initial $\bar{\mathbf{L}}_n$, respectively.

## 4. Speed Adjustment on the Given DBLs

### 4.1. State Definition

*4.1.1. State Variable*

A state variable is a series of features of a bus line corresponding to a given critical time point. To choose a proper control action at a critical time point, we need to pick out some key features of the bus line to form the state variable. Let $t_m$ be the critical time point in question. The system state of a bus line should have the potential to reproduce the trajectories of buses with the help of control actions. In view of the above consideration, we will choose three features of any bus $b$ to constitute the state variable.

The first feature is the current dwelling stop for a stopped bus or the next stop to be reached for a cruising bus. This feature will be called the **target stop** for the bus in question. $e_b^m$ is used to denote the target stop of bus $b$ with respect to the critical time point $t_m$. The second feature is the **arrival time** for a bus to arrive at its target stop. $t_b^{\mathbb{A},m}$ is used to denote the arrival time of bus $b$ at stop $e_b^m$ with respect to the critical time point $t_m$. The third feature is the time interval from the current time to the departure time regarding the target stop $e_b^m$. This will be called **the time interval for a bus to be activated.** $t_b^{\mathbb{D},m}$ is used to denote the time interval for bus $b$ to be activated at the critical time point $t_m$. Group all the above three features for all buses at the critical time point $t_m$ into a vector $s_m$. $s_m$ is the state variable with respect to the critical time point $t_m$. Use Figure 1 as an example to demonstrate the above conceptions. Assume that the current time and the current critical time point $t_m$ with respect to $b_1$ are $t_1$ and $t_3$, respectively. We can see that $e_{b_1}^m$, $t_{b_1}^{\mathbb{A},m}$ and $t_{b_1}^{\mathbb{D},m}$ are $e_1$, $t_2$ and $t_3 - t_1$, respectively.

*4.1.2. How to Estimate the Expected Dwell Time*

To carry out the multi-CTP look-ahead operation to be presented in the next subsection, we need to estimate the dwell time for a bus at its target stop before the bus arrives at the stop. Assume that we are considering bus $b$ and stop $e$ here. The passengers generation rate with respect to the stop $e$ is denoted by $r_e$. The average boarding time for a passenger is $t_{board}$. Assume that the latest bus left stop $e$ at $t_D$. The estimated arrival time of bus $b$ at stop $e$ is $t_{b,e}^A$. So we can use $\max\{(t_{b,e}^A - t_D)r_e t^{board}(1 + r_e t^{board}), t_{b,e}^{alight}\}$ as the estimate of the dwell time of bus $b$ at stop $e$. $(t_{b,e}^A - t_D)r_e t^{board}$ is the boarding time required for passengers who arrive at $e$ before bus $b$ arrives. $(t_{b,e}^A - t_D)(r_e t^{board})^2$ is the estimate of the boarding time required for passengers who arrive at stop $e$ during the dwelling period of bus $b$. $t_{b,e}^{alight}$ is the estimated alighting time of passengers in bus $b$ at stop $e$. To obtain the estimate of $t_{b,e}^{alight}$, we need to know the average alighting time of a passenger and the total number of passenger on bus $b$ who will alight at $e$. This will be realized in the numerical experiment section by a detailed simulation system.

### 4.2. Cost of action

When a control action is determined, its impact on the stability of bus line needs to be measured so as to judge the quality of the action. With the DCH defined by Eq. (1), we can define the cost of action $a$ as follows

$$c(a) = \sum_{b \in B}(h_b(t) - H(t))^2, \qquad (6)$$

where the headway $h_b(t)$ is corresponding to the critical time point $t$ when the control action $a$ is determined. $H(t)$ stands for the DCH corresponding to time $t$. The cost $c(a)$ has the same changing trend as $\sigma_H(t)$ defined in eq. (2).

Since $a$ stands for the regulating speed on the dedicated bus lanes, its impact on $h_b(t)$ comes from the changing travel times of the dedicated bus lanes. But there is one situation which needs to be further clarified here. When there is no speed adjustment at a critical time point (e.g. the action set is {0}), Eq. (6) will still be used to calculate the influence of some action earlier determined. In this case, we should view $c(a)$ as a general conception used to evaluate the stability of the bus line at time $t$.

### 4.3. Speed Adjustment by Multi-CTP Look-ahead

To choose an action, that is a regulating speed in our study, by looking ahead several CTPs is to roll forward several CTPs of the simulation system of the bus line to see which action will bring us the best result. To roll the simulation system forward, we need to use the expected values for many quantities, such as the expected delay at an intersection and the dwell time at a bus stop.

Assume that the current time is $\vec{t}$, the current system state is $s_m$ corresponding to CTP $t_m$ and the number of successive CTPs to be looked ahead is $N$. The series of CTPs to be looked ahead is $\{t_1, t_2, \ldots, t_i, \ldots t_N\}$. The explanations of some new notations are as follows. $e_b^i$ stands for the target stop for bus $b$ to be activated when the current looking ahead CTP is $t_i$. $t_b^{D,i}$ is the time remaining for bus $b$ to be activated at the looking ahead CTP $t_i$. $t_e^{la,i}$ is the modified latest arrival time corresponding to bus stop $e$ when the current CTP is $t_i$. Note that we use the time point $\vec{t}$ as the origin to define the modified latest arrival time. To distinguish the expected value from the actual value, the overbar will be used to indicate the expected values. So $\overline{t}_g^b$ stands for the expected travel time of the bus line segment $g$ for bus $b$. $\overline{t}_e^b$ is the expected dwell time of bus $b$ at stop $e$. $\Delta \overline{t}_g^b(a)$ is used to denote the changed travel time of the bus line segment $g$ for bus $b$ due to the control action $a$.

The process of searching a proper action $a^*$ with $N$-CTP look-ahead is as follows. Note that if the control action set corresponding to the current CTP is $\{0\}$, the process will be overlooked. We use $e \oplus 1$ to stand for the nearest downstream stop of bus stop $e$ in the following process.

**Step 1**: Assign a large positive value to $\tilde{c}_i, \forall i \in \{1, 2, \cdots, N\}$.

**Step 2**: Choose a bus $b^1 \in \arg\min_{d \in B}\{t_d^{\mathbb{D},m}\}$ to be active at the first level.

**Step 3:** Carry out the following multiple nested "for" loops.

For every action $a^1 \in A_{e_{b^1}^m}$, execute the following operations (start the first level):

$$e_{b^1}^1 := e_{b^1}^m \oplus 1 \text{ and } e_b^1 := e_b^m, \forall b \neq b^1;$$

$$t_{b^1}^{D,1} := t_{b^1}^{\mathbb{D},m} + \Delta \bar{t}_{g^1}^{b^1}(a^1) + \bar{t}_{g^1}^{b^1} + \bar{t}_{e_{b^1}^1}^{b^1} \text{ and } t_b^{D,1} := t_b^{\mathbb{D},m}, \forall b \neq b^1;$$

$$t_{e_{b^1}^1}^{la,1} := t_{b^1}^{\mathbb{D},m} + \Delta \bar{t}_{g^1}^{b^1}(a^1) + \bar{t}_{g^1}^{b^1} \text{ and } t_e^{la,1} := t_e^{\mathbb{A},m} - \vec{t}, \quad \forall e \neq e_{b^1}^{w,1}.$$

Note that the two ends of the bus line segment $g^1$ are $e_{b^1}^m$ and $e_{b^1}^1$.

Calculate the cost of action $c(a^1)$ by eq. (6).

Choose a bus $b^2 \in \arg\min_{d \in B}\{t_d^{D,1}\}$ to be activated at the second level.

For every $a^2 \in A_{e_{b^2}^1}$, execute the following operations (start the second level):

$$e_{b^2}^2 := e_{b^2}^1 \oplus 1 \text{ and } e_b^2 := e_b^1, \forall b \neq b^2;$$

$$t_{b^2}^{D,2} := t_{b^2}^{D,1} + \Delta \bar{t}_{g^2}(a^2) + \bar{t}_{g^2} + \bar{t}_{e_{b^2}^2}^{b^2} \text{ and } t_b^{D,2} := t_b^{D,1}, \forall b \neq b^2;$$

$$t_{e_{b^2}^2}^{la,2} := t_{b^2}^{D,1} + \Delta \bar{t}_{g^2}(a^2) + \bar{t}_{g^2} \text{ and } t_e^{la,2} := t_e^{la,1}, \quad \forall e \neq e_{b^2}^2.$$

Calculate the cost of action $c(a^2)$ by eq. (6).

…

Choose a bus $b^N \in \arg\min_{d \in B}\{t_d^{D,N-1}\}$ to be activated at the level $N$.

For every $a^N \in A_{e_{b^N}^{N-1}}$, execute the following operations (start the level $N$):

$$e_{b^N}^N := e_{b^N}^{N-1} \oplus 1 \text{ and } e_b^N := e_b^{N-1}, \forall b \neq b^N;$$

$$t_{b^N}^{D,N} := t_{b^N}^{D,N-1} + \Delta \bar{t}_{g^N}(a^N) + \bar{t}_{g^N} + \bar{t}_{e_{b^N}^N}^{b^N} \text{ and } t_b^{D,N} := t_b^{D,N-1}, \forall b \neq b^N;$$

$$t^{la,N}_{e^N_{b^N}} := t^{D,N-1}_{b^N} + \Delta \bar{t}_{g^N}(a^N) + \bar{t}_{g^N} \quad \text{and} \quad t^{la,N}_e := t^{la,N-1}_e, \forall e \neq e^N_{b^N}.$$

Calculate the cost of action $c(a^N)$ by eq. (6);

Let $\tilde{c}_N := \min\{\tilde{c}_N, c(a^N)\}$.

**End** the "for" of the level $N$.

Let $\tilde{c}_{N-1} := \min\{\tilde{c}_{N-1}, c(a^{N-1}) + \gamma \tilde{c}_N\}$.

…

Let $\tilde{c}_2 := \min\{\tilde{c}_2, c(a^2) + \gamma \tilde{c}_3\}$.

**End** the "for" of the level 2.

Let $\tilde{c}_1 := \min\{\tilde{c}_1, c(a^1) + \gamma \tilde{c}_2\}$ and $a^* := \arg\min_{a^1}\{\tilde{c}_1, c(a^1) + \gamma \tilde{c}_2\}$

**End** the "for" of the level 1.

**Step 4**: Output the optimal action $a^*$ with respect to the current state $s_m$.

In the above process, $\gamma \in (0,1]$ is a given rate. $\gamma$ is used to discount the expected cost.

Note that during the above multi-CTP look-ahead process, the corresponding action set at some critical time points is $\{0\}$. In this case, the corresponding 'for' loops will degenerate into simply carrying out the rolling operation with no speed adjustment.

## 5. Numerical Experiments

### 5.1. Information about the Bus Line to Be Tested

To test our method, a detailed circular bus line with 36 stops and 11 cruising buses is constructed. The length of the bus line is 21.35km. Table 1 shows the passenger capacity, the initial stops where buses are dwelling at or will reach first at the beginning of the simulation, and the Time Remaining for a bus to be Activated at the Beginning of every simulation (TRAB). Note that the time unit is second (s) in the following.

Table 1
The basic data about the buses.

| No. of Bus | 1 | 2 | 3 | 4 | 5 | 6 | 7 | 8 | 9 | 10 | 11 |
|---|---|---|---|---|---|---|---|---|---|---|---|
| Capacity | 72 | 70 | 80 | 60 | 72 | 60 | 72 | 80 | 60 | 72 | 70 |
| Initial stop | 1 | 4 | 8 | 11 | 15 | 18 | 21 | 25 | 28 | 31 | 34 |
| TRAB (s) | 20 | 0 | 40 | 30 | 50 | 10 | 30 | 36 | 24 | 18 | 26 |

The basic data of the bus line segments are shown in Table 2. 'BLS' and 'RS' in Table 2 stand for Bus Line Segment and Road Segment, respectively. The inclusion relations between BLS and RS and the actual lengths of RS are listed. We use red color and boldface to indicate the possible alternative locations where DBLs may be installed.

Table 2

The relations between bus line segment and road segment.

| BLS | RS | Length (m) | BLS | RS | Length (m) | BLS | RS | Length (m) |
|---|---|---|---|---|---|---|---|---|
| 1 | 1,2 | 200,400 | 13 | 18,19 | 300,320 | **25** | **36** | **530** |
| **2** | **3** | **500** | 14 | 20 | 500 | 26 | 37 | 560 |
| **3** | **4** | **600** | 15 | 21 | 450 | 27 | 38 | 600 |
| 4 | 5,6 | 260,350 | 16 | 22,23,24 | 200,250,100 | 28 | 39,40 | 300,500 |
| **5** | **7** | **530** | **17** | **25** | **570** | **29** | **41** | **600** |
| 6 | 8 | 560 | 18 | 26 | 610 | 30 | 42,43 | 300,350 |
| 7 | 9 | 600 | 19 | 27 | 600 | 31 | 44,45 | 200,400 |
| 8 | 10,11 | 300,500 | **20** | **28,29** | **300,350** | 32 | 46 | 500 |
| 9 | 12 | 600 | **21** | **30,31** | **200,400** | **33** | **47** | **600** |
| 10 | 13,14 | 300,350 | 22 | 32 | 500 | **34** | **48,49** | **260,350** |
| **11** | **15** | **600** | 23 | 33 | 600 | 35 | 50 | 530 |
| 12 | 16,17 | 300,400 | 24 | 34,35 | 260,350 | 36 | 51 | 560 |

Note. 'BLS' and 'RS' stand for bus line segment and road segment, respectively.

We assume that the average cruising speed of buses in a common road segment without DBL is 35(km/hr), but the average speed of buses in a dedicated bus lane is 50(km/hr). These average speeds will be used in the multi-CTP look-ahead process. To mimic the actual running of the bus line, we need to consider the uncertainty of travel times. Assume that the length of the road segment $d$ is $l_d$ (km). With the average cruising speed $\bar{v}$ in the road segment, the expected travel time $\bar{t_d}$ of $d$ is $l_d / \bar{v}$. A sample travel time of $d$ can be generated by $t_d = \bar{t_d} + \delta_d$ where $\delta_d$ is a normal random variable with zero mean and variance $\sigma_d^2$. The obtained sample travel time will be used as the actual travel time in our study. To simplify the computation, we assume that $\sigma_d = 5l_d$ holds for a common road segment and $\sigma_d = 2l_d$ for a DBL.

There are 15 intersections along the bus line. We assume that the signal control schemes at these intersections are all pre-timed two phases regarding the bus-coming approaches. The red phase $t_i^{red}$, the green phase $t_i^{green}$, and the remaining time of the initial phase $t_i^{or}$ are listed in three rows. The row labeled by "Initial phase" specifies the initial phases of all the intersections. In this row, 1 and 2 indicate the red and green phases, respectively. The last row points out the corresponding bus line segment which includes the given intersection.

Table 3

The basic data of intersections.

| No. of Intersection | 1 | 2 | 3 | 4 | 5 | 6 | 7 | 8 | 9 | 10 | 11 | 12 | 13 | 14 | 15 |
|---|---|---|---|---|---|---|---|---|---|---|---|---|---|---|---|
| $t_i^{red}$ (s) | 40 | 40 | 40 | 30 | 30 | 40 | 40 | 30 | 30 | 40 | 40 | 40 | 30 | 40 | 40 |
| $t_i^{green}$ (s) | 50 | 30 | 35 | 45 | 30 | 30 | 45 | 35 | 45 | 50 | 30 | 35 | 45 | 50 | 30 |
| $t_i^{or}$ (s) | 20 | 20 | 10 | 20 | 20 | 20 | 30 | 20 | 20 | 10 | 20 | 10 | 20 | 20 | 20 |
| Initial phase | 2 | 1 | 1 | 2 | 2 | 1 | 2 | 2 | 2 | 2 | 1 | 1 | 2 | 2 | 1 |
| Bus line segment | 1 | 4 | 8 | 10 | 12 | 13 | 16 | 16 | 20 | 21 | 24 | 28 | 30 | 31 | 34 |

Passengers will be generated by a Monte Carlo simulation according to the given arrival rates. The corresponding relationships between stops and arrival rates are shown in Table 4. We also assume there are two types of passengers specified by their average boarding and alighting times. The average boarding and alighting times of the first type are 4s and 1s, respectively. The average boarding and alighting times of the second type are 2s and 0.5s, respectively. Assume that 10% of all the passengers are grouped into the first type.

Table 4

Bus stops with given arrival rate.

| $r_e$ | Related Stops |
|---|---|
| 1 | 4,6,8,9,12,13,14,16,17,20,21,24,26,28,30,32,34,36 |
| 2 | 1,2,3,5,10,22,23,25,31,35 |
| 3 | 7,11,18,19,27,29,33 |
| 4 | 15 |

Note. The unit of $r_e$ is passenger per minute in this table.

Table 5

The probabilities for the following downstream stops to be chosen as destination.

| No. | 1 | 2 | 3 | 4 | 5 | 6 | 7 | 8 | 9 | 10 | 11 | 12 | 13 |
|---|---|---|---|---|---|---|---|---|---|---|---|---|---|
| Series 1 | 0.0135 | 0.027 | 0.0541 | 0.0811 | 0.1081 | 0.1351 | 0.1351 | 0.1216 | 0.1216 | 0.0811 | 0.0541 | 0.0405 | 0.0270 |
| Series 2 | 0.0345 | 0.0862 | 0.1207 | 0.1552 | 0.1724 | 0.1552 | 0.1207 | 0.0862 | 0.0517 | 0.0172 | / | / | / |

When a passenger is generated, the corresponding destination stop can be determined by a given probability distribution related to the following alternative stops. Two series of probabilities are given in Table 5. Series 1 and 2 have 13 and 10 elements, respectively. The sum of the values of the elements in any one series is 1. The $n$ th element of the series used by the current stop is the probability of the $n$ th downstream stop to be chosen as a destination stop. We assume that series 1 is used by bus stops including 1, 7, 10, 13, 20, 21, 27, 30, 31 and series 2 is used by the other remaining stops.

At the end of the observation, we need to evaluate the service level of the bus line from the viewpoint of passengers. We will use $t_P^W$, $t_P^R$ and $t_P^{Tr}$ to denote the average waiting time, the average riding time and the average travel time over all passengers $P$ who have finished their bus trip during the observation period. $\sigma_P^W$, $\sigma_P^R$ and $\sigma_P^{Tr}$ denote the standard deviations of $t_P^W$, $t_P^R$ and $t_P^{Tr}$ over the set $P$, respectively. $n_P$ is the number of passengers in $P$.

In view of the applicability in practice, we assume that the action set of the speed adjustment is {-10, -5, 0, 5, 10} with measurement unit km/hr. This set will be used throughout the numerical experiments. The discount rate $\gamma$ is 0.5. To obtain the expected value of $\mathbb{E}[\bar{c}_H(x, A^\pi(s,x))]$, we will carry out 50 times of the simulation and then use the average value of $\bar{c}_H(x, A^\pi(s,x))$ over these 50 times simulations as the expected value. The other quantities to be shown in the subsequent tables are also the average values obtained in a similar way as above. The observation time period is 4 hours.

**5.2. Stabilizing the bus line with the given DBLs**

*5.2.1. An Initial Crude impression*

# Position of figure 3

Fig. 3 The trajectories of buses with no control.

# Position of figure 4

Fig. 4 The trajectories of buses with the speed adjustment control.

The bus line in question has a strong trend of bus bunching under the non-control scenario. A typical result of buses' trajectories in this scenario is drawn in Figure 3. We can see the bus bunching becomes more and more serious at the end of the observation period of 4 hours.

Figure 4 presents a typical result of buses trajectories obtained with the speed adjustment control on the alternative locations indicated in Table 2. The number of CTPs to be looked ahead in the control scheme is 3 in this case. Contrary to the performance in Figure 3, there is no bus bunching appearing in Figure 4. The buses cruise in a relatively smooth way throughout the observation period.

*5.2.2. The Influence of the Number of CTPs to Be Looked Ahead*

In this subsection, we will investigate the influence of the specified number of CTPs to be looked ahead in the control scheme. We still consider 11 bus line segments to install DBLs highlighted in Tables 2. $\sigma_{\bar{c}}$ is the standard deviation of $\bar{c}_H$ on the set $\bar{T}$ such that

$$\sigma_{\bar{c}} = \sqrt{[\sum_{t\in\bar{T}}(\sigma_H(t)-\bar{c}_H)^2]/(n_{\bar{T}}-1)}, \tag{7}$$

where $n_{\bar{T}}$ is the number of CTPs in $\bar{T}$. $\sigma_{\bar{c}}$, i.e. the standard deviation of $\bar{c}_H$, shows the reliability of $\bar{c}_H$ as an expected estimate.

Since the regulating speed can be either positive or negative. To properly evaluate the strength of the speed adjustment, we will use the absolute value of regulating speed to form our measurement values. $a_\Sigma$ is the sum of all the absolute values of the regulating speeds adopted during our control period. $\bar{a}$ is the average value of all the absolute values of the regulating speeds. $\sigma_{\bar{a}}$ is the deviation of the absolute values of the regulating speeds to $\bar{a}$. $n_M$ is the number of the regulating speeds adopted during our control period.

With the above defined notions, we give out in Tables 6 and 7 the results with respect to different numbers of CTPs to be looked ahead in the speed adjustment. From the data in Table 6, we can see the performance of 3CTPL outperforms the others with respect to the stability index $\bar{c}_H$. All the results of the 5 types of speed adjustment schemes outperform the non-control scheme. The bus bunching is removed from the bus line with the speed adjustment as indicated in the column labeled by 'Bunch'. The absolute values of the average regulating speeds for these control schemes are similar to each other. The value of $\bar{c}_H$ decreases with the increased number of CTPs at the beginning, but the trend changes to increase later with the continually increased number of CTPs to be looked ahead.

Table 6
The stability indices and control actions with respect to different strategies.

| Methods | $\bar{c}_H$ | $\sigma_{\bar{c}}$ | $n_{\bar{T}}$ | $a_\Sigma$ | $\bar{a}$ | $\sigma_{\bar{a}}$ | $n_M$ | Bunch |
|---|---|---|---|---|---|---|---|---|
| No Control | 496.93 | 388.24 | 2083 | / | / | / | / | Yes |
| 1CTPL | 37.8 | 9.92 | 2112 | 2995 | 4.63 | 8.35 | 647 | No |
| 2CTPL | 35.8 | 8.47 | 2097 | 2925 | 4.56 | 8.38 | 641 | No |
| 3CTPL | 33.4 | 6.44 | 2111 | 3005 | 4.65 | 8.30 | 646 | No |
| 4CTPL | 34.7 | 8.70 | 2119 | 2940 | 4.54 | 8.23 | 648 | No |
| 5CTPL | 37.0 | 8.75 | 2110 | 3140 | 4.87 | 8.61 | 645 | No |

Note. 'nCTPL' stands for the strategy of looking ahead $n$ Critical time points.

Table 7
Waiting and riding times of passengers with respect to different strategies.

| Methods | $n_P$ | $t_P^W$ | $\sigma_P^W$ | $t_P^R$ | $\sigma_P^R$ | $t_P^{Tr}$ | $\sigma_P^{Tr}$ |
|---|---|---|---|---|---|---|---|
| No Control | 14089 | 353.1 | 309.5 | 427.7 | 193.4 | 780.8 | 368.0 |
| 1CTPL | 14538 | 134.6 | 73.0 | 414.8 | 189.0 | 549.4 | 202.5 |
| 2CTPL | 14674 | 133.8 | 72.5 | 428.4 | 193.3 | 562.3 | 206.8 |
| 3CTPL | 14784 | 131.7 | 71.8 | 415.2 | 188.7 | 546.9 | 203.0 |
| 4CTPL | 14449 | 132.7 | 72.3 | 418.5 | 187.8 | 551.2 | 202.9 |
| 5CTPL | 14572 | 131.6 | 72.7 | 417.0 | 188.3 | 548.5 | 202.6 |

The data in Table 7 show that all the speed adjustment strategies outperform the non-control scheme with respect to the waiting and riding times of passengers. But among these speed adjustment strategies, no one dominates others in all aspects with respect to the time attributes of passengers. For example, 5CTPL has minimal waiting time 131.6s, but 1CTPL has minimal riding time 414.8s.

The above observation shows that there generally exists an optimal number of CTPs used in the multi-CTP look-ahead. The influence of the different number of CTPs on the service level is generally negligible.

*5.2.3. The Influence of the Number of DBLs to Be Installed*

The branch and bound algorithm presented in Section 3 assumes that with the reduced number of locations for installing DBLs, the capacity of the speed adjustment strategy in the DBLs will decrease. So as a result of the above assumption, the instability trend will increase with the reduced number of locations. To test this assumption, we will consider 5 sets of locations where DBLs are installed. The specific sets are given in Table 8 where 'BLS' stands for the bus line segment that is the location with installed DBLs.

Table 8

The sets of bus line segments where DBLs are installed.

| Sets of BLSs | The corresponding Bus line segments |
|---|---|
| 11BLS | 2,3,5,11,17,20,21,25,29,33,34 |
| 9BLS | 2,3,5,11,17,20,21,25,29 |
| 7BLS | 2,3,5,11,17,20,21 |
| 5BLS | 2,3,5,11, 17 |
| 3BLS | 2,3,5 |

Table 9

The stability indices and control actions regarding different DBLs' locations.

| Methods | $\bar{c}_H$ | $\sigma_{\bar{c}}$ | $n_{\bar{T}}$ | $a_\Sigma$ | $\bar{a}$ | $\sigma_{\bar{a}}$ | $n_M$ | Bunch |
|---|---|---|---|---|---|---|---|---|
| 11BLS | 35.8 | 8.47 | 2097 | 2925 | 4.56 | 8.38 | 641 | No |
| 9BLS | 36.8 | 8.93 | 2116 | 2590 | 4.90 | 8.52 | 529 | No |
| 7BLS | 44.3 | 11.70 | 2093 | 2125 | 5.21 | 8.97 | 408 | No |
| 5BLS | 52.3 | 13.62 | 2084 | 1845 | 6.34 | 10.01 | 291 | No |
| 3BLS | 153.7 | 167.79 | 2125 | 1505 | 5.97 | 9.80 | 252 | Yes |

Note. 'nCTPL' stands for the strategy of looking ahead n critical time points.

Table 10

Waiting and riding times of passengers regarding different DBLs' locations.

| Methods | $n_P$ | $t_P^W$ | $\sigma_P^W$ | $t_P^R$ | $\sigma_P^R$ | $t_P^{Tr}$ | $\sigma_P^{Tr}$ |
|---|---|---|---|---|---|---|---|
| 11BLS | 14674 | 133.8 | 72.5 | 428.4 | 193.3 | 562.3 | 206.8 |
| 9BLS | 14517 | 131.5 | 73.0 | 415.3 | 186.4 | 546.9 | 201.0 |
| 7BLS | 14687 | 136.8 | 77.0 | 417.9 | 191.4 | 554.8 | 207.7 |
| 5BLS | 14324 | 140.5 | 81.3 | 419.4 | 191.1 | 559.9 | 209.0 |
| 3BLS | 14279 | 165.3 | 107.7 | 437.8 | 269.2 | 603.1 | 286.2 |

Tables 9 and 10 present the computation results with respect to the sets in Table 8. Note that we will look ahead two CTPs in the following computations. Three observations are worth noting here. Firstly, the changing trends of the data of $\bar{c}_H$ and $\sigma_{\bar{c}}$ in Table 9 are corresponding to our assumption mentioned above. This guarantees the soundness of our designed branch and bound algorithm. Secondly, when the size of the set of the chosen locations becomes small and the locations, that are the elements of the set, are relatively near to each other, the speed adjustment strategy may be out of work. The set of locations {2,3,5} makes the strategy powerless to remove the bus bunching as indicated in the column labeled by 'Bunch' in Table 9. Thirdly, the service levels with respect to the waiting and riding times of passengers show a changing trend. At the beginning, the service level increases and then changes to decrease gradually. Note that this changing trend does not correspond to the changing trend of the stability index. The reason behind the last observation needs some further study in the future.

**5.3. Determine the Optimal Locations of DBLs**

*5.3.1. Demonstration of the Searching Process*

Table 11
The $F_i^1$ and $F_i^2$ of alternative locations.

| No. of BLS | 2 | 5 | 11 | 17 | 20 | 21 | 25 | 29 | 33 | 34 |
|---|---|---|---|---|---|---|---|---|---|---|
| $F_i^1$ | 2.5 | 2.65 | 3.0 | 2.85 | 3.25 | 3.0 | 2.65 | 3.0 | 3.0 | 3.05 |
| $F_i^2$ | 12.28 | 10.95 | 8.23 | 9.33 | 6.59 | 8.23 | 10.95 | 8.23 | 8.23 | 7.88 |

In this subsection, we will use a simple example to show the process of the branch and bound algorithm. The values of $F_i^1$ and $F_i^2$ with respect to the alternative locations to be considered are presented in Table 11. To make the demonstration clear, we consider a relatively small set of locations {2, 5, 17, 20, 25} at the beginning. Set the upper limit $\hat{F}^1$ to 10 and $\hat{F}^2$ to a very large positive number. The searching process of the optimal locations in set {2, 5, 17, 20, 25} is shown in Figure 5.

# Position of figure 5

Fig. 5 The searching tree with the root node of {2, 5, 17, 20, 25}.

In the searching tree of Figure 5, the serial numbers of the chosen feasible locations for a node to install DBLs are given in the big bracket '{}' beside the node. The number tailed by a capital letter '*F*' or '*T*' beside a node is the objective value of the node. Here '*T*' stands for that the constraint to the influence of installing DBLs on traffic flow is satisfied; '*F*' means the constraint is not satisfied. Any directional line with arrow connects a parent node to one of its descendant nodes. The two numbers labeled beside a directional line are the orders with which the line is traversed during the searching. For example, the numbers 2(3) beside the line connecting nodes {5, 17, 20, 25} to {17, 20, 25} means that during the searching, the line has been traversed from {5, 17, 20, 25} to {17, 20, 25} and then from {17, 20, 25} to {5, 17, 20, 25} at the 2nd and 3rd steps in a searching order, respectively. Five alternative locations have 32 combinations. From Figure 5, we know 13 combinations have been checked.

*5.3.2. The Influence of Constraints on Results*

We will use the 10 locations given in Table 11 as the root node in the subsection. For a set with 10 elements, we know it has 1024 different combination. In other words, if we check all the combinations, we need to generate and check 1024 nodes.

Table 12

The optimal locations and objectives with respect to different constraints.

| Scenario | $\hat{F}_1 >= F^1 x^*$ | $\hat{F}_2 >= F^2 x^*$ | Objective | Optimal Locations | Nodes Generated | Feasible Nodes |
|---|---|---|---|---|---|---|
| 1 | **25>=23.6** | **70>=69.29** | 30.45 | 5,11,20,21,25,29,33,34 | 36 | 14 |
| 2 | **20>=20** | 80>65.48 | 33.47 | 2,5,11,17,21,29,33 | 168 | 56 |
| 3 | 20>16.65 | **60>=59.97** | 37.38 | 2,5,11,17,25,29 | 221 | 79 |
| 4 | **16>=14.55** | 60>=44.62 | 38.85 | 5,17,21,33,34 | 347 | 69 |
| 5 | 18>=14.15 | **50>=47.69** | 38.06 | 5,17,21,25,29 | 312 | 43 |
| 6 | **14>=14.0** | **50>=49.02** | 39.07 | 2,11,17,25,29 | 560 | 52 |
| 7 | **13>=11.55** | 45>=36.39 | 40.60 | 5,17,29,34 | 618 | 63 |

By changing the upper limits of constraints, we investigate 7 scenarios as shown in Table 12. The column labeled by 'Nodes Generated' lists the number of nodes generated during the searching. The column labeled by 'Feasible Nodes' lists the number of nodes whose set of the chosen location satisfies the constraints. From the data in Table 12, we can see that with the decrease of two upper limits, the objective values and the number of nodes generated during the searching increase. The number of the corresponding optimal locations also decreases with the decrease of the upper limits. If only one constraint is in consideration and no more nodes can be added in the set of the chosen locations, the constraint can be called 'tight'; or else, the constraint can be called 'relaxed'. According to the above definition, we highlight the tight constraints with red color and boldface in Table 12.

By carefully observing the results obtained until now, we find out that with given constraints, the set of the optimal locations should generally have the most elements and the elements should evenly disperse along the bus line in most cases. An obvious example is that the result regarding set {2, 3, 5} given in Table 9 is contrary to the result regarding set {5, 17, 25} shown in Figure 5. This observation can be explained as follows. As we have checked in subsection 5.2.3, the more DBLs are installed, the stronger the capacity of the speed adjustment strategy will be. So an optimal set of the chosen locations should include the most possible number of locations. Here we overlook the feature of the location, such as the length of the corresponding bus line segment, by assuming that these features are similar with respect to different locations. Why should the feasible locations used to install DBLs disperse evenly along the bus line in most cases? If these locations are concentrated together, the range of the influence of the speed adjustment will be relatively small along the bus line comparing to the deployment of evenly dispersed locations.

**5.4. Computation Time**

The program used in this section is written in Java 1.8.0_45 and runs in NetBeans IDE 8.0.2. The computer used is a laptop with a processor of Intel® Core i3-3120M CPU @2.50GHz and installed memory (RAM) of 4.00GB (2.32GB usable).

Using the feasible locations given in Table 11, the times to simulate the bus line operation in a control period of 4 hours for strategies 1CTPL, 2CTPL, 3CTPL, 4CTPL, and 5CTPL are 0.156s, 0.25s, 0.735s, 3.121s, and 11.187s, respectively. The times to supply a regulating speed at the critical time points regarding strategies 1CTPL, 2CTPL, 3CTPL, and 4CTPL are less than 0.001s. The times to supply a regulating speed for strategies 5CTPL is about 0.012s. These results are very promising for the implementation of our methods in practice.

Denote the time to simulate the bus line operation in a control period of 4 hours by $t_S$ and the time to supply a regulating speed by $t_K$. Obviously, $t_S$ can be approximated by $n_M t_K$ where $n_M$ is the number of the actual speed adjustments during the observation time period. The value of $n_M$ is determined by the length of the observation period and the number of locations used to install DBLs. Use $n_{Tree}$ to denote the number of nodes generated during a searching process. Use $n_{EV}$ to denote the number of runs of simulations specified to calculate the expected objective value. We can use $n_{Tree} n_{EV} t_S$ to estimate the total computation time required to identify the optimal locations from the given alternative locations. Since $t_S$ is relatively small as shown in the preceding paragraph. We can expect the final total computation time will be acceptable in most of the cases.

## 6. Conclusion

In this paper, we have addressed two problems. One is to determine the optimal locations to install DBLs with the given funding limit and the maximal allowable influence on the other road traffic. We have formulated a non-linear mixed integer stochastic mathematical model and design a specific branch and bound algorithm to solve this problem. The other is to adjust the speeds of buses in these DBLs to improve the operational stability of the bus line. A multiple Critical Time Point look-ahead process is proposed to determine the optimal regulating speed. The solution of the second problem supplies the estimate of the influence of DBLs on the other traffic. This estimate is required to solve the first problem.

The numerical experiments have verified the proposed branch and bound algorithm to our location selection model and the speed adjustment strategy. The results of the numerical analyses show that the new branch and bound algorithm can considerably reduce the number of combinations of the feasible locations to be checked. The speed adjustment strategy can stabilize a high-frequency bus line which has a strong trend to bunching.

Three important observations from the numerical analyses are summarized as follows. Firstly, the number of CTPs to be looked ahead should not be too big. Due to the various stochastic factors, the estimates of various expected values will become inaccurate with the increased number of CTPs. According to the simulation experiments, 3 CTPs to be looked ahead is better than the other choices. In practice, the specified number of CTPs should be determined carefully by trial and error. Secondly, we find out that the more DBLs are installed along a bus line, the more powerful the speed adjustment in these DBLs will be with respect to stabilizing an unstable bus line. This observation confirms the common sense. But we should note that to guarantee the above conclusion, a proper method of adjusting the bus speeds in DBLs is required. The third observation is about two features of the resulted optimal deployment of DBLs. One feature is that under the given constraints, the optimal deployment of the chosen locations to install DBLs should include the feasible locations as many as possible. The other feature is that the chosen locations should be evenly distributed along the bus line. The above observation can be used directly by practitioners to deploy DBLs in real life.

There are several interesting directions which can be further pursued in the future. In view of the influences among buses servicing different bus lines and cruising in the common DBLs, considering the deployment of DBLs in a bus line network will improve the effectiveness of the implementation of related strategies. Combining the speed control in DBLs with transit signal priority at intersections can be a natural extension since they are usually considered together by public transit system managers in practice.

**Acknowledgment**

This research was supported in part by National Natural Science Foundation of China(71601118, 71801153, 71871144），the Natural Science Foundation of Shanghai(18ZR1426200) and the Key Climbing Project of USST (No. SK17PA02).

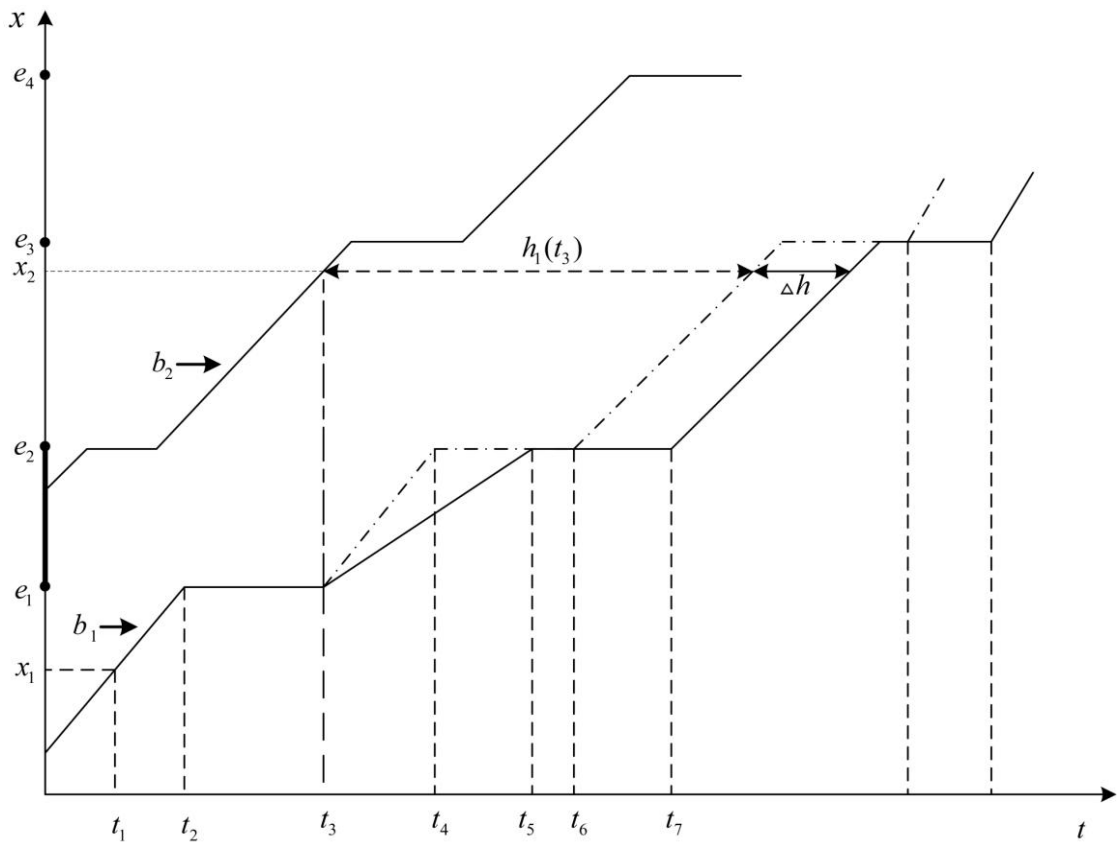

Fig. 1

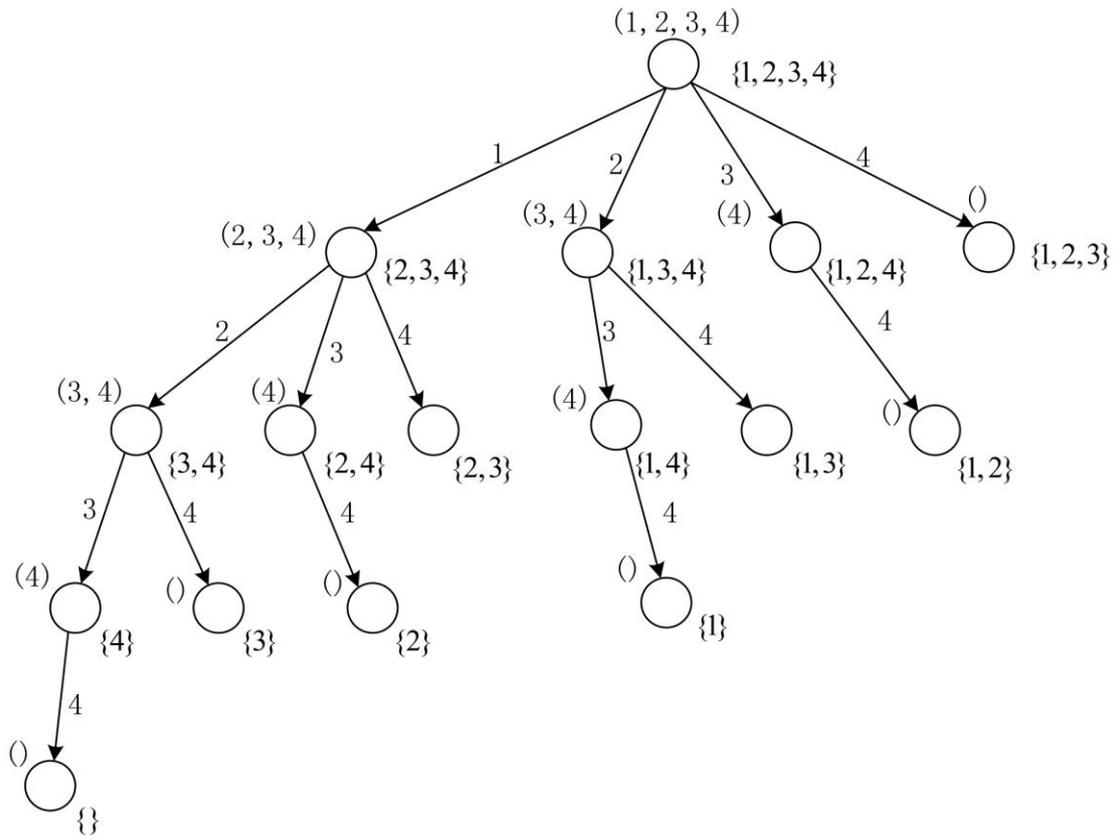

Fig.2

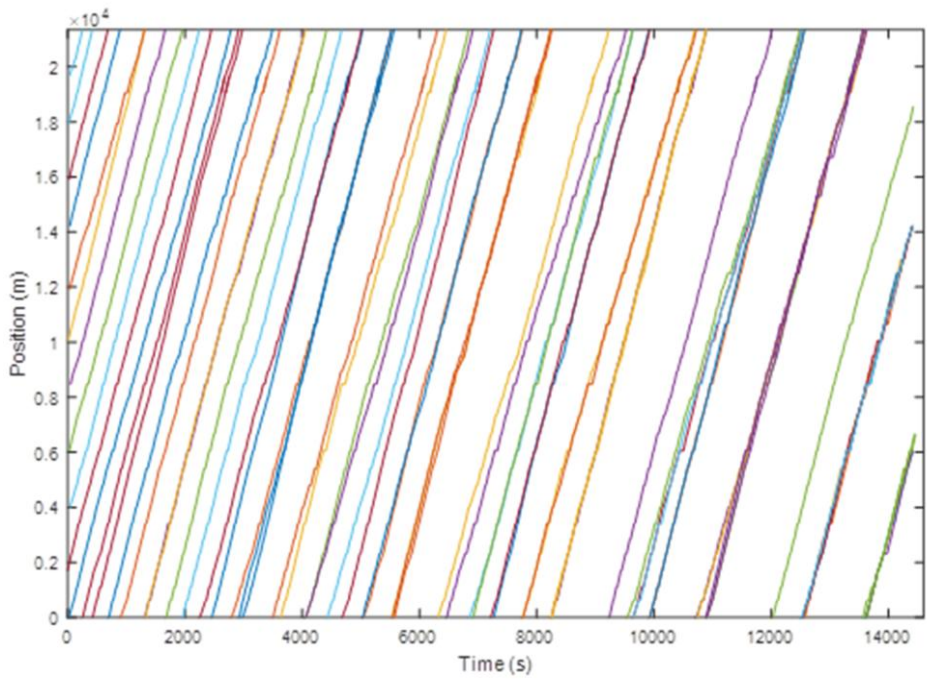

Fig.3

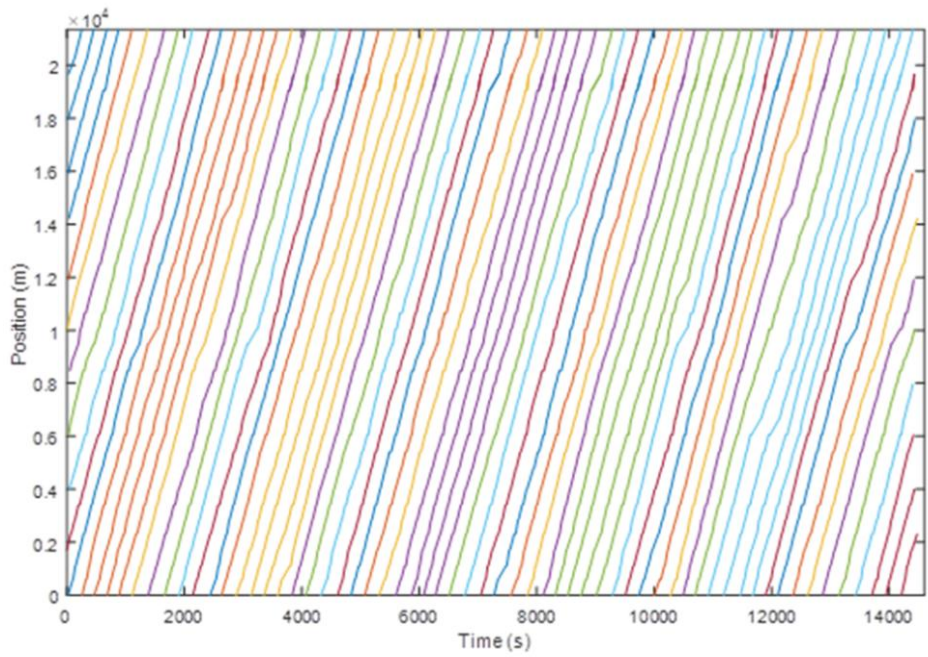

Fig. 4